\documentclass[11pt]{article}
\usepackage{amssymb, amsmath, amsfonts}
\newtheorem{theorem}{Theorem}[section]
\newtheorem{lemma}[theorem]{Lemma}

\newtheorem{proposition}[theorem]{Proposition}
\newtheorem{corollary}[theorem]{Corollary}
\newtheorem{conjecture}[theorem]{Conjecture}

\newcommand{\comment}[1]{}

\title{A near neighbour continuum percolation model}
\author{Alexis Gillett and Misja Nuyens \\  Vrije Universiteit Amsterdam}

\begin{document}
\maketitle

\begin{abstract}

We introduce a continuum percolation model defined on the points of a $d$-dimensional homogeneous Poisson process. Each Poisson point is connected to all points within its connection range, which depends on the distances to the other Poisson points. We show that the new model exhibits a phase transition, and obtain results about the critical values in low and high dimensions.

\end{abstract}


\section{Introduction}

A continuum percolation model consists of a point process and a rule for connecting the points. The first mention of these types of models appears to date from 1961 \cite{gilbert}. Typically, the point process is  a homogeneous Poisson process $X$ on $\mathbb{R}^d$ with some given density $\lambda$, and the points are connected in one of two ways. In a {\em random connection} model, points  are connected to each other by undirected edges determined by some probability measure. In a {\em Boolean} model, a ball is placed around each point in $X$, and its radius is generated by some probability measure. Two points are considered connected if their balls overlap, and clusters are formed in the obvious way.
An overview of these types of models can be found in  \cite{roy}.

In their 1996 paper on continuum percolation \cite{hag}, Meester and H\"aggstr\"om looked at some  models for which the density of the underlying Poisson process is irrelevant to the percolation probability. In these models, altering the density is equivalent to rescaling the model. In this paper, we generalise this work. Let us begin by reviewing the original models.

The {\em hard sphere Boolean} models are those Boolean models that assign probability $1$ to the configurations where none of the balls overlap, although they are allowed to touch. One such example is the dynamic lily-pond model, where the configuration of balls is generated by growing balls at each point linearly in time. As soon as a ball touches another ball, it stops growing. The existence of this model is non-trivial since there is a.s.\ no first contact. It has  been shown that this model exists and that there is no percolation in any dimension \cite{hag}.

The {\em nearest neighbour}
 model is an example of a random connection model. The model is defined by the following connection rule: connect each $x \in X$ to its $k$ nearest neighbours. Letting $U$ denote the event that an infinite component is formed, we define
$$k_c(d)=\min \{k \ge 1: \mathbb{P}(U)>0 \}.$$
By ergodicity, we have that $k_c(d)=\min \{k \ge 1: \mathbb{P}(U)=1 \}$. It has  been proven for this model that $k_c(1)=\infty$ and that $2 \le k_c(d) <\infty$, for $d \ge 2$, see \cite{hag}. Furthermore, it is known that $k_c(d) = 2$ for large $d$. It is also believed, but unproven, that $k_c(2)=3$ and $k_c(d)=2$ for $d \ge 3$. We shall denote the nearest neighbour model by NN$(d,k)$, where $d$ is the dimension and $k \in \mathbb{N}$ is the parameter.

In the next section we introduce a generalised nearest neighbour model, of which the NN($d,k)$ model is a special case.
The main questions that we shall study in this new model concern the percolation probability, and in particular the existence of phase transitions as the dimension of the Poisson process (or any other parameter for that matter) is varied. An important tool at our disposal is renormalisation. Renormalisation means discretising the space into boxes to form a lattice and then looking at the connections between neighbouring boxes. This then permits a comparison with site percolation. There are a number of examples of this kind of approach \cite{sar, mass}. Some of our proofs will be based on the methods used in \cite{hag}.

\section{A generalised nearest neighbour model}\label{sec:2}

Let $X$ be a homogeneous Poisson process on $\mathbb{R}^d$.
Consider a point $x \in X$ and let $d_i(x)$ denote the Euclidean distance from $x$ to its $i$th nearest neighbour. The NN$(d,k)$ model can now be defined by the connection rule ``$x$ is connected to all points within distance $d_k(x)$ of $x$.'' We want to generalise this model by not considering the $d_i(x)$ themselves, individually, but a function of all these distances. Since we want the model to be density-invariant, it is clear that this function should be linear in that $f(\{d_i(x)\}_{i \ge 1})=\sum_{i=1}^{\infty}f_i(d_i(x))$, with the $\{f_i\}_{i \ge 1}$ being linear maps from $\mathbb{R}^+$ to $\mathbb{R}^+$, i.e., $f_i(x)=\alpha_i x$.

With this motivation, we define a new model by connecting $x$ by an edge to
all points $y$ with $|x-y| \le r(x)$, where $r(x)$ is given by
\[r(x) = \sum_{i=1}^{\infty} \alpha_i d_i(x),\]
for some vector  $\underline{\alpha}=(\alpha_1,\alpha_2,\dots)$ with $\alpha_i\geq 0$ for all $i$.
We call this model the {\em generalised nearest neighbour} model, and denote this model by GN($d,\underline{\alpha})$. Since we shall begin our analysis by looking at simpler versions of this model, we denote by GN$_k(d,\alpha)$ the model when only the $k$th component of $\underline{\alpha}$ is non-zero, and has  value $\alpha$. Observe that GN$_k(d,1)$ = NN$(d,k)$.

Looking at the GN$_1(d,\alpha)$ model, we define a critical value for $\alpha$ as before:
$$\alpha_c(d)=\inf \{\alpha \ge 0: \mathbb{P}(U)>0 \}.$$
Since the critical value for NN$(d,1)$ is known to be $\infty$, we have that $\alpha_c(d)\geq 1$ for all $d$. Note that an alternative approach would be to fix $\alpha$ and look for the critical value of $d$. One could also consider a corresponding Boolean model, where on each point $x$ of the Poisson process a sphere  of radius $\alpha d_1(x)$ is placed. This model is related to previously considered dependent continuum percolation models: for example, for $\alpha \ge 1$, the  Boolean model dominates all hard sphere models, since in  a hard sphere model  the ball at $x$ must necessarily have radius less than or equal to $d_1(x)$.

For the GN$_k(d,\alpha)$ model, we denote the corresponding critical value by $\alpha_c^{(k)}(d)$. Hence under this notation, $\alpha_c(d)=\alpha_c^{(1)}(d)$. Sometimes, however, we will vary only $k$ and look for a phase transition for GN$_k(d,\alpha)$ with $\alpha$ and $d$ fixed. Note the following property when $\alpha<1$: for all points $x \in X$, there are at most $k-1$ points within distance $r(x)=\alpha d_k(x)<d_k(x)$. This leads to the result that the GN$_k(d,\alpha)$ model is dominated by the GN$_{k-1}(d,1)$ model for all $\alpha<1$. In the  GN$(d,\underline{\alpha})$ model there is no obvious definition of a critical value as we are dealing with a (possibly infinite dimensional) vector. However, one number that will be of interest is $|\underline{\alpha}|=\sum_{i=1}^{\infty}\alpha_i$. Note that percolation is  trivial in all dimensions if $|\underline{\alpha}|=\infty$. This case is therefore ignored, which restricts us to $\underline{\alpha}$ such that $\alpha_i \to 0 \text{ as } i \to \infty$.

In the formation of large clusters, there are two phenomena that are competing. First, to form large clusters, points in the cluster should be connected to (many) other points, and to do this, points should be close to each other. On the other hand, points should have a large connection range, so points should be far away from each other. As we shall see in the next section, this makes the model already far from trivial in one dimension. This is in contrast to  the NN($d,k$) models, for which it is very  easy to show that there is no percolation in one dimension.

Although the generalised nearest neighbour model is a random connection model, it is possible to define a corresponding Boolean model. In the Boolean version, two points are connected if their connection ranges overlap. The proofs of a number of our results can be easily modified to apply to this Boolean model, for example Theorem \ref{banana} below. Typically any numerical bounds on the critical value, such as Corollary \ref{number}, will be halved since the connection ranges now combine to make connections.

This Boolean model is of particular interest when we set $\alpha_1=1$ and $\alpha_i=0$ for all $i>1$. In this case, the Boolean version of the generalised nearest neighbour model dominates all homogeneous Poisson hard-sphere models. It is an interesting open problem to determine in what dimensions percolating hard sphere models exist. The existence of such models has recently been proven for $\mathbb{R}^d$ when $d \ge 45$ \cite{cotar} and the existence is also known on some spaces other than $\mathbb{R}^d$ \cite{jonasson}. On $\mathbb{R}$, it is known that no such models exist, so that the problem is still open for $2 \le d \le 44$. Showing that $\alpha_c(d)>1$ (for some $d$) for the Boolean version of the generalised nearest neighbour model would prove that there is no percolation in $d$-dimensional hard-sphere models. However, numerical simulations for $d=2$ suggest that this approach will not be successful.

The paper is organised as follows.  We first study the GN model in one dimension, in Section \ref{subs1}. Then in Sections \ref{subs2} we consider the GN$_1(d,\alpha)$ model, where only the distance to the nearest neighbour is important.
In Section \ref{subs4} we study the GN$_k(d,\alpha)$ model for large $d$ and in  Section \ref{subs5} we obtain some results for the general GN($d,\underline{\alpha})$. The paper concludes with the statement of some open problems.

\section{The model in dimension 1}\label{subs1}

For the nearest neighbour model NN(1,$k$) there is no percolation in one dimension for any $k$. This result is a straightforward consequence of the following two features. First, the Poisson process on $\mathbb{R}$ contains arbitrarily large gaps. Second, if there is no edge between two neighbouring points, then these points belong to separate clusters. For the GN$_k(1,\alpha)$ model though, this second property does not hold and the proof that there is no percolation in one dimension, given below, is non-trivial.

We first introduce some notation.
Let $(x,x+m)$ be called a {\em gap of length $m$} if $X \cap (x,x+m) = \emptyset$. The term $m$-gap denotes a gap of length greater than $m$. We say that there is a {\em bridge} over a gap if two points on different sides of the gap have an edge between them; we say  that $x$ {\em bridges} a gap if the connection range of $x$, $r(x)$, spans that gap. Furthermore, a gap $(x,x+m)$ is bridged from the right if there is a point $y\geq x+m$ such that $y-r(y)<x$. The properties of homogeneous Poisson processes imply that for any point $x \in X$, there are infinitely many $m$-gaps in both the positive (to the right) and negative (to the left) directions for all $m<\infty$.

Let $p(m)$ denote the probability that an $m$-gap is not bridged from the right.

\begin{lemma}\label{lem:pb}  Consider the GN$_k(1,\alpha)$ model, and let $\beta> (\alpha\vee 1)$. Then $p(\beta^2)>0$.
\end{lemma}
{\bf Proof: }
It will be convenient to use  a second  distance function, $\bar{d}_k(x)$, where $\bar{d}_k(x)$ is the distance from $x$ to its $k$th nearest neighbour to the right.
Furthermore,  call $x \in X$ a {\em $\beta$-point} if $\bar{d}_k(x)>\beta$. Observe that if $x$ is not a $\beta$-point, then $r(x) \le \alpha \bar{d}_k(x) < \alpha \beta < \beta^2$. Hence, when looking for bridges from the right over $\beta^2$-gaps, it is sufficient to only consider $\beta$ points.

Consider a $\beta^2$-gap, and the first point to the right of it. Without loss of generality, call this point the origin, $0$. Let $0=X_0 < X_1 < X_2 < \dots$ denote all points to the right of $0$ and let $0 \le Y_0 < Y_1 < Y_2 < \dots$ denote the $\beta$-points to the right of $0$. A sufficient condition for the $\beta^2$-gap to be unbridged from the right is that $Y_i-\alpha \bar{d}_k(Y_i) > - \beta^2$ for all $i \ge 0$. From the definition of $\beta$-points it follows that $Y_{i+k}-Y_i>\beta$, and hence $Y_i \geq  \lfloor \frac{i}{k}\rfloor \beta$.

Since the events $\{\alpha \bar{d}_k(Y_i)\leq t\}$ are positively correlated for all  $i$ and $t$, we have
\begin{align}
p(\beta^2) & \geq  P(\alpha \bar{d}_k(Y_ i)\leq Y_i+\beta^2 \ \textrm{for all} \ i)\nonumber\\ &\geq
 P(\alpha \bar{d}_k(Y_ i)\leq \lfloor\frac{i}{k}\rfloor \beta +\beta^2 \ \textrm{for all} \ i)\nonumber\\
& \geq \prod_{i=1}^{\infty} P(\alpha \bar{d}_k(Y_ i)\leq \lfloor\frac{i}{k}\rfloor \beta +\beta^2 ).
\label{eq.een}\end{align}
Let $\Gamma$ have  a gamma$(k,1)$ distribution, and note that $\bar d_k(x)\stackrel{d}{=}\Gamma$. Since $Y_i$ is a $\beta$-point, and $\beta>(\alpha\vee 1)$, we have the following for all $i>k\beta$:
\begin{align}
 P(\alpha \bar{d}_k(Y_ i)\leq \lfloor\frac{i}{k}\rfloor\beta +\beta^2 ) &
 \geq P(\bar{d}_k(Y_ i)\leq \lfloor\frac{i}{k}\rfloor  +\beta )
\geq P(\bar{d}_k(Y_ i)\leq \frac{i}{k} )
\nonumber\\
& \geq  P( \Gamma   \leq  \frac{i}{k} \ |\ \Gamma>\beta)  =1-P(\Gamma   \geq  \frac{i}{k} \ |\ \Gamma>\beta).
\label{eq.twee}\end{align}
We now use  that for $n\to\infty$,
\begin{equation}\label{eq.drie} P(\Gamma \geq n|\Gamma>\beta)=\frac{\sum_{i=1}^{k-1}e^{-n}\frac{n^i}{i!}}
{\sum_{i=1}^{k-1}e^{-\beta}\frac{\beta^i}{i!}}\sim c(\beta,k) n^{k-1}e^{-n},\end{equation}
where $1/c(\beta,k)=\sum_{i=1}^{k-1}e^{-\beta}(k-1)!\beta^i/i!$.
Combining (\ref{eq.een}), (\ref{eq.twee}) and (\ref{eq.drie}),
we find that $p(\beta^2) > \prod_{i=1}^{\infty} (1-a_i)$
with $a_i\sim c'(\alpha, \beta,k) i^{k-1}e^{-i}$ as $i\to\infty$ and $c'(\alpha, \beta,k)>0$,
and $a_i>0$ for all $i$.
Since $a_i \to 0$ exponentially fast, we conclude that $p(\beta^2)>0$. \hfill $\Box$

\begin{theorem} \label{prop1}
In the GN$_K(1,\alpha)$ model,
we have  $\alpha_c^{(k)}(1)=\infty$.
\end{theorem}
{\bf Proof:}
 We prove the theorem by showing that there exist  a.s.\ infinitely many unbridged $\beta^2$-gaps in both directions.
By symmetry,  the probability that a $\beta^2$-gap is bridged from the left is $p(\beta^2)$ as well.
Now note that for a gap of given length, between two Poisson points $x$ and $y$ say, the event that it is unbridged from the right is positively correlated with the event that it is unbridged from the left. Indeed,  the absence of a  bridge from the right makes the  distances of points to the right of $y$ to their $k$th nearest neighbour stochastically smaller, and as a consequence, the same holds for  the distances from the points to the left of the gap to their $k$th nearest neighbour.
Therefore, the probability that a $\beta^2$-gap is unbridged is at least $p(\beta)^2$, which is strictly positive by Lemma \ref{lem:pb}.

This would complete the proof if $\beta^2$-gaps were bridged independently. Since this is not the case, we describe a scanning procedure to demonstrate that the desired result holds. Assume without loss of generality that there is a Poisson point at the origin.
Start at the origin and scan to the right until you find a $\beta^2$-gap.
This will be unbridged with probability $p(\beta^2)$.
If the $\beta^2$-gap is not bridged, then this makes the neighbouring $\beta^2$-gaps more likely to be unbridged, by the same reasoning as above. We then continue scanning to the right of the gap.

If there is a bridge over the $\beta^2$-gap, originating at the Poisson point $x$ say, then we jump ahead to $x+r(x)$,  walk $k$ points to the right, and continue scanning. Since the dependence between the connection ranges of points does not go further than their $k$ nearest neighbours, this removes the dependence on the possible bridging of the previous gap(s).

By the  positive correlation  (previous gap was unbridged) or independence (previous gap was bridged), all gaps encountered by the scanning procedure  are unbridged with probability at least $p(\beta)^2$. Hence, we can recover an infinite sequence of $\beta^2$-gaps such that each is unbridged with probability at least $p(\beta^2)$, and the events that gaps are unbridged are either independent, or positively correlated.
Hence, there are a.s.\ infinitely many unbridged $\beta^2$-gaps to the right.
The proof is completed by noting that the scanning procedure to the left works exactly the same.  \hfill $\Box$\\

We now show that if the $\alpha_i$ are large enough, then the GN$(1,\underline{\alpha})$ model does percolate. In fact, we show that for these $\alpha_i$, the model is {\em fully connected.}

\begin{proposition}
In the GN$_k(1,\alpha)$ model, let $x$ be a point of the Poisson process on $\mathbb{R}$. If $\sum_{i=1}^\infty i\alpha_i=\infty$, then $r(x)=\infty$ a.s.
As a consequence, the GN$(1,\underline{\alpha})$ model is fully connected.
\label{warm}
\end{proposition}
{\bf Proof: }
First note that $(d_1(x), d_2(x)-d_1(x), \ldots)\stackrel{d}{=}(U_1,  U_2, \ldots)$, where $U_1, U_2, \ldots$ are i.i.d.\  exponential random variables with parameter 2.
Hence, defining $d_0(x)=0$, we may write
\begin{align} r(x)&=\sum_{i=1}^\infty \alpha_id_i(x) =\sum_{i=1}^\infty \alpha_i\sum_{j=1}^i d_{j}(x)-d_{j-1}(x)\stackrel{d}{=} \sum_{i=1}^\infty \alpha_i\sum_{j=1}^i U_j\nonumber\\
& =U_1\sum_{i=1}^{\infty} \alpha_i + U_2\sum_{i=2}^{\infty} \alpha_i+\cdots\nonumber
=\beta_1U_1+\beta_2U_2+\cdots,\label{alphabeta}\end{align}
where $\beta_k=\sum_{i=k}^\infty \alpha_i$. Before continuing, note that
$\sum_i\beta_i=\sum_ii\alpha_i$.
Denote the Laplace transform of a random variable $X$ by $\phi_X$.
Since the $U_i$ are independent, the Laplace transform of $V=\sum_i \beta_iU_i$ satisfies
\[ \phi_{V}(s)=\prod_{i=1}^\infty \phi_{U_i}(\beta_is)=
\prod_{i=1}^\infty \frac{1}{1+\beta_is/2}=
\prod_{i=1}^\infty \Big(1-\frac{\beta_is}{2+\beta_is}\Big).\]
Since the $U_k$ are independent and $\{V=\infty\}$ is a tail event, its probability is either 0 or 1, by a Zero-One Law. Hence, $V=\infty$ if and only if  $\phi_{V}(s)=0$ for all $s> 0$.
But $\phi_{V}(s)=0$ if and only if $\sum_i\beta_i=\infty$.
Hence, if  $\sum_i\beta_i=\infty$, then   $V=\infty$ a.s., and therefore $r(x)=\infty$ a.s.~\hfill $\Box$\\

There is an obvious potential extension of this result to higher dimensions.
This result is given in Section \ref{subs5}.

Observe that $\mathbb{E}(r(x))=\mathbb{E}(U_1)\sum_{i=1}^{\infty} \alpha_i + \mathbb{E}(U_2)\sum_{i=2}^{\infty} \alpha_i+\cdots = \frac{1}{2}\sum_{i=1}^\infty i\alpha_i$. Therefore, if the conditions of Proposition \ref{warm} are not satisfied, the expected connection range is finite. Furthermore, it can be shown that the second moment of $r(x)$ can be bounded from above by a polynomial of its first moment. Therefore, in this case, the variance is also finite. For related one-dimensional independent percolation models, we have the result that there is no percolation if $\mathbb{E}(r(x))<\infty$ \cite{hall}. Thus, it seems reasonable to conjecture that the converse of Proposition \ref{warm} is true.
\begin{conjecture}\label{1d}
If $\sum_{i=1}^\infty i\alpha_i<\infty$, then all clusters in the GN$(1,\underline{\alpha})$ model are finite almost surely.
\end{conjecture}
Currently, we have only the following weaker Theorem.

\begin{theorem}\label{t.expo}
If $\alpha_i=\gamma^i$ for some $0<\gamma\leq 1/2$, then the GN$(1,\alpha)$ model does not percolate.
\end{theorem}
{\bf Proof:}
Since the events that gaps are not bridged are positively correlated, it suffices to show that  there exists an $m$ such that the probability that an $m$-gap is not bridged is strictly positive.
Since the events that a gap is not bridged from the left and from the right are positively correlated, it suffices to show that the probability that a gap is not bridged from the right is strictly positive.

To do that, denote the points of the Poisson process on $[0,\infty)$ by $X_0\leq X_1\leq \cdots$. Then for all $i$ and $k$ we have
\[ d_i(X_k)\leq d_i(X_0)+X_k-X_0.\]
Hence, we can calculate
\begin{align*} r(X_k)& = \sum_i \gamma^i d_i(X_k)\leq \sum_i \gamma^i [d_i(X_0) +X_k-X_0]\\
&= r(X_0)+\frac{\gamma}{1-\gamma} [X_k-X_0]\leq r(X_0)+X_k-X_0.
\end{align*}
Hence, for all $k$,
\[ X_0-r(X_0)\leq X_k-r(X_k).\]
So, if $X_0$ does not bridge the gap to its left, then neither do $X_1, X_2,\ldots$. To complete the proof, it suffices to show that there exists an $m$ such that $P(r(X_0)< m)>0$. Letting $\bar d_i$ denote the distance to the $i$th nearest neighbour on the right, we calculate,
\begin{align*} P(r(X_0)<m)& =P\left( \sum_i \gamma^i d_i(X_0)<m \right) \geq
P\left(\sum_i \gamma^i \bar d_i(X_0)<m\right)\\ &=P\left(\sum_i \gamma^i (U_1+\cdots +U_i)<m\right)=
P\left(\sum_k \frac{\gamma^k U_k}{1-\gamma}<m\right),\end{align*}
where the $U_i$ are i.i.d.\ exponentially distributed with parameter 1.
Hence, by the Markov inequality,
\[ P(r(X_0)<m)\geq 1- P\left(\sum_{k=1}^\infty \frac{\gamma^k U_k}{1-\gamma}\geq m\right)
\geq 1- \frac{1}{m}E\sum_{k=1}^\infty \frac{\gamma^k U_k}{1-\gamma}=1-\frac{\gamma}{m(1-\gamma)^2}.\]
We conclude that $P(r(X_0)<m)>0$ for $m>\gamma/(1-\gamma)^2$, which completes the proof.\hfill $\Box$\\

Obviously, Theorem \ref{t.expo} also holds for all $\underline{\alpha}$ with $\alpha_i\leq \gamma^i$ for all $i$ and some $0< \gamma\leq 1/2$.  Furthermore,
readers familiar with long-range percolation may also see superficial similarities with the one-dimensional homogeneous case where $q(n)<1$ denotes the probability of being connected to a vertex at distance $n$. If $\sum n q(n) < \infty$, there is no percolation \cite{grim}.

\section{The GN$_1(d,\alpha)$ model  for $d\geq 2$}\label{subs2}

Having completed our treatment of the special case $d=1$, we continue by looking at $d \ge 2$, starting with the simplest model: GN$_1(d,\alpha)$. The next result shows that there exists a non-trivial critical value for all $d>1$.

\begin{figure}[htb]
\begin{center}
\begin{picture}(150,150)(0,0)

\put(0, 0){\line(1,0){150}}
\put(0, 0){\line(0,1){150}}
\put(150, 0){\line(0,1){150}}
\put(0, 150){\line(1,0){150}}

\put(0, 60){\line(1,0){60}}
\put(60, 0){\line(0,1){60}}

\put(20, 20){\line(1,0){20}}
\put(20, 20){\line(0,1){20}}
\put(20, 40){\line(1,0){20}}
\put(40, 20){\line(0,1){20}}

\put(35, 35){\circle*{3}}
\put(55, -10){$3$}
\put(145, -10){$3n$}

\end{picture}
\end{center}
\caption{A banana box inside a good box of side length $3n$.}
\label{figgie}
\end{figure}
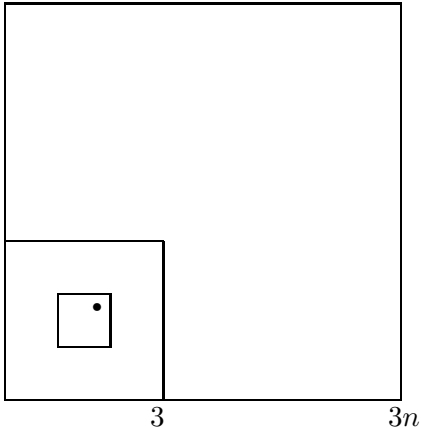

\begin{theorem} For $d\geq 2$, $\alpha_c(d)<\infty$.
\label{banana}
\end{theorem}
{\bf Proof}
We first consider $d=2$, with  the density of the Poisson process   equal to $1$.
A $3\times 3$ box is called a {\em banana box} if the $1\times 1$  box in its centre contains exactly one point, and the rest of  box is empty,
see  Figure \ref{figgie}. The probability that a  $3\times 3$ box is  a  banana box is $e^{-1}e^{-8}$.
A $3n\times 3n$ box is called {\em good} if it contains at least one banana box. The probability that a  $3n\times 3n$ box is good is at least $1-(1-e^{-9})^{n^2}.$

Now choose $n$ so large that the probability that a $3n\times 3n$ box is good is larger than $p_c^{site}$. Consider two neighbouring good boxes. By construction, in a good box there is a point whose nearest neighbour is at distance at least 1. On the other hand, two  points in two neighbouring $3n\times 3n$ boxes are at most $3n\sqrt{5}$ away from each other.  Then for all $\alpha>n\sqrt{45}$, any two neighbouring good boxes will be connected to each other. Since the probability that a box is good is larger than $p_c^{site}$, there is a.s.\ an infinite cluster of good boxes. Hence, $\alpha_c(2)\leq n\sqrt{45}< \infty$. For  $d\geq 3$, the proof is similar, and is therefore omitted. $\hfill \Box$\\

We now use the same idea to find a  not so sharp upper bound for $\alpha_c(2)$.
In the  proof of Theorem \ref{banana} we used unit boxes purely for convenience.  To alter the scale of our boxes, we now introduce an extra parameter, $\delta$. This $\delta$ allows us to optimise the choice of $n$, so that we can   minimise our upper bound for $\alpha_c(2)$. A $3\delta\times 3\delta$ box is called a $\delta$-banana box if the inner box of side length $\delta$ contains one point, and the rest of the box is empty. This happens with probability $\delta^2 e^{-9\delta^2}$. An $3\delta n \times 3\delta n$ box is called $\delta$-good if it contains at least one $\delta$-banana box. So, if we want the probability that a $3\delta n \times 3\delta n$ box is $\delta$-good to be larger than $p_c=p_c^{site}$, it suffices to have

\begin{equation}\label{bnd} 1-(1-\delta^2 e^{-9\delta^2})^{n^2}\geq p_c.\end{equation}

Since our bound for $\alpha_c(2)$ will be $n\sqrt{45}$, we want $n$ in (\ref{bnd}) to be as small as possible. To do that, we maximise the  LHS of (\ref{bnd}). Analytically, we readily find that $\delta=1/3.$
There is also a probabilistic way to see this. For a $3\delta\times 3\delta$ box to be a $\delta$-banana box, it must contain only one point, and this point must be in the central square. Given that there is only one point, the probability that it is in the central square is $1/9$ for all $\delta$. If we wish to maximise the probability that a $3\delta\times 3\delta$ box is a $\delta$-banana box,  we merely need to maximise the probability that there is only one point in the box. This leads to $\delta=1/3$, so that the smallest $n$ that satisfies \ref{bnd} is given by
\[ \tilde{n} = \Big\lceil \sqrt{\frac{\log (1-p_c)}{\log (1- \frac{1}{9e})}} \Big\rceil.\]
Monte Carlo simulations suggest that the critical value for site percolation on the square lattice $\mathbb{Z}^d$ is approximately $0.59$ \cite{ziff}. Since the exact value is unknown, in the following corollary  we use the best available rigorous upper bound: $p_c\leq  0.679492$ \cite{wierman}.
\begin{corollary} We have
$ \alpha_c(2)\leq \tilde{n}\sqrt{45} < 41.$
\label{number}
\end{corollary}

We now move to the asymptotical behaviour of $\alpha_c(d)$ for $d\to\infty$.

\begin{theorem}\label{muchwork} We have
$\alpha_c(d) \to 1$  as  $d \to \infty$.
\label{biggie}
\end{theorem}

The proof of Theorem \ref{muchwork} builds on the work of Meester and H\"aggstr\"om \cite{hag} and follows a similar approach to Section 3 of their paper. They varied $k$ rather than $\alpha$ and looked at the critical value of $k$ for the GN$_k(d,1)$ model. In their paper, it was proven that $k=2$ is the critical value when $d$ is large. It should be noted that the result was not originally presented in this form and has been expressed here in the terminology of this paper. The technique they use is showing that there is percolation in the GN$_2(d,1)$ model for all sufficiently large $d$. Combining this with an earlier result that the GN$_1(d,1)$ model doesn't percolate for any $d$ yields the result. Here, we derive Theorem \ref{muchwork} by showing that for any $\alpha > 1$  the GN$_1(d,\alpha)$ model percolates for $d$ large enough.

The proof of this result requires a number of steps. Before giving the formal proof, we give an outline of the proof and some preliminary results. Our approach will be to construct a point process that with non-zero probability gives a subset of the cluster at the origin of GN$_1(d,\alpha)$ that is infinite.  The existence of such a process would implies that the GN$_1(d,\alpha)$ model percolates. The point process is built up by a sequence of steps. A step is said to be {\em successful} if it produces a (finite) subset of the cluster at the origin of GN$_1(d,\alpha)$ and initiates two subsequent steps. The idea is that by taking the dimension very large, we can uniformly bound from below the probability that a step is successful. By making this lower bound sufficiently large, we then show that an infinite sequence of steps occurs with positive probability.

The {\em out-cluster} of a point $x$ in GN$_1(d,\alpha)$ is the set of points defined by the following iterative procedure, initiated by the set $\{x\}$. Given that we have a set $\{x_1,x_2,\dots,x_n\}$, add to this set all Poisson points $y$ such that $d(x_i,y) < \alpha d_1(x_i)$ for some $i$. This is then repeated with the new set of points. The procedure either continues forever, in which case the out-cluster is infinite, or stops when there are no more points to add. By construction, all points in the out-cluster at $x$ must belong to the cluster at $x$. We shall prove percolation for the GN$_1(d,\alpha )$ model by showing that (a subset of) the out-cluster at the origin is infinite with positive probability.

We attempt to construct a subset of the out-cluster at the origin via a point process. This point process is built up algorithmically using $d$-dimensional spatial branching processes. Initially, the space ($\mathbb{R}^d$) is empty and the spatial branching processes place points in this space to create the point process. Each spatial branching process (SBP) is run for $n$ generations and is thought of as a step in the algorithm. A step of the algorithm is declared {\em successful} if two things happen. First, the union of the SBP with all previously successful steps must form a subset of the out-cluster at the origin of GN$_1(d,\alpha)$. This means that for all points $x$ in the point process, the offspring of that point must be within distance $\alpha d_1(x)$ of $x$. Also, the collection of points must obey the law of a homogeneous Poisson process. Second, there must be two points in the $n$-th generation of the SBP that satisfy a certain location condition as defined below. These two points are used to initiate two more steps of the algorithm. The point process is then defined as the union of all the successful steps of the algorithm.

The SBP lives in $\mathbb{R}^d$, but one can consider its projection onto $\mathbb{R}^2$ via the linear map
\begin{equation}\label{weekend} L(x_1,\dots,x_d)=\sqrt{d}(x_1,x_2).\end{equation}
This permits a comparison with oriented site percolation on the lattice ${\cal L}=\{(i,j)\in \mathbb{Z}^2 : i \ge 0, |j| \le i, (i+j)/2 \in \mathbb{Z}\}$, with oriented edges from $(i,j)$ to $(i+1,j\pm 1)$. Each site $(i,j)$ of ${\cal L}$ corresponds to a square $S_{i,j}=[M(i-\frac{1}{2}),M(i+\frac{1}{2})] \times [M(j-\frac{1}{2}),M(j+\frac{1}{2})]$, where $M$ will be chosen later. The second condition for a step to be successful (mentioned above) is now that the projected SBP starting in box $S_{i,j}$ has points in its $n$-th generation in $S_{i+1,j-1}$ and $S_{i+1,j+1}$. This allows us to start a pair of new SBPs from these points with projected origins in $S_{i+1,j-1}$ and $S_{i+1,j+1}$ respectively. In this way, a successful step in the algorithm corresponds to an open site in ${\cal L}$.
Let $p_c<1$ be the critical value of oriented site percolation on ${\cal L}$.
If for all sufficiently high dimensions we can bound uniformly (i.e., irrespective of what has happened in previous steps of the algorithm) from below the probability of a successful step of the algorithm by some $p>p_c$, , then we have shown that the GN$_1(d,\alpha)$ model percolates. Before proving this result, we define the spatial branching process and give some properties.\\

\noindent {\em The spatial branching process. (SBP)} Let $S_r(x)$ denote the hyper-sphere of radius $r$ centred at $x$, write $S_r=S_r(0)$, and let $|S_r|$ denote the volume of $S_r$. The spatial branching process with origin $0$ is defined as follows. We start with $Z_0=\{0\}$.  Given $Z_n$, the  offspring of each $y \in Z_n$ is generated by the following procedure. Let $\delta_1>0$. An independent homogeneous Poisson point process $X^{(y)}$ on $S_{1+\delta_1}(y)$ is generated, with density $\lambda(d)$ such that the expected number of points contained in $S_1(x)$ is $1$. Then a ball is grown around $y$ until either the ball has radius $1+\delta_1$ or $c_2$ points of $X^{(y)}$ have been encountered, for a certain $c_2\in\mathbb{N}$. The offspring of $y$ are then the points of  $X^{(y)}$ that
are contained in this ball.

In general, the SBP does not generate a homogeneous Poisson process. Consider two points $x$ and $y$ such that $d(x,y)<1+\delta_1$. It is clear that the overlap $S_{1+\delta_1}(x) \cap S_{1+\delta_1}(y)$ is non-empty, and that considering the offspring of both $x$ and $y$ together, the density of points is doubled on this overlap. However, if we condition on $y$ having no offspring on this overlap, then the joint collection of offspring does form a homogeneous Poisson process $S_{1+\delta_1}(x) \cup S_{1+\delta_1}(y)$. Thus if the overlap is small, then with high probability, namely when $y$ has no offspring on this overlap, we can consider the union of the offspring to be a homogeneous Poisson process on $S_{1+\delta_1}(x) \cup S_{1+\delta_1}(y)$. This reasoning will be crucial when comparing the SBP to the out-cluster at the origin of the GN$_1(d,\alpha)$ model. The following standard result shows that in high dimensions, this overlap is negligible.

\begin{lemma}\label{noverlap}
If $x_1$ and $x_2$ are such that $d(x_1,x_2)>0.9$ and $r_1,r_2 \in (0.9,1.1)$, then
$$\frac{|S_{r_1}(x_1) \cap S_{r_2}(x_2)|}{|S_{r_1}(x_1)|} \to 0$$
as $d \to \infty$.
\end{lemma}

Observe that for all $\epsilon > 0$, $|S_{1+\epsilon}|/|S_1| \to \infty \text{ as } d \to \infty$. The offspring distribution for each individual in the SBP is {\em distributed like} $Y \wedge c_2$, where $Y$ is  Poisson distributed with parameter $|S_{1+\delta_1}|/|S_1|$. Hence, for any $0<c_1<c_2$, we can chose $\delta_1=\delta_1(d)$ such that for $d$ sufficiently large, $\mathbb{E}(Y \wedge c_2)=c_1$  and that $\delta_1(d) \to 0 \text{ as } d \to \infty$. Note that $c_1$ and $c_2$ do not vary with the dimension. Thus with this choice of $\delta_1(d)$, the offspring of a point $x$ converge weakly (as $d\to\infty$) to a set of points uniformly distributed on the surface of the unit sphere centred at $x$. If $c_1>1$, then the SBP is supercritical and the probability of extinction can be made arbitrarily small by taking $c_1$ sufficiently large (independent of the dimension).

The following lemma comes from \cite{pen}, and used the map $L$ defined in (\ref{weekend}).

\begin{lemma}
Suppose $U=(U_1,\dots,U_d)$ is uniform on the surface of $S_1$. Then, as $d \to \infty$, the two-dimensional random vector $L(U)$ converges in distribution to the bivariate normal distribution $N(0,I)$ with zero mean and as covariance matrix the identity matrix $I$.
\label{tod2}
\end{lemma}

The mapping $L$ can be used to map our SBP onto $\mathbb{R}^2$. Let SBP$^*$ denote the limit as $d \to \infty$ of the mapped SBP. Since $L$ is continuous, this limit is the same as taking the mapping the limit of the SBP. Thus, by the previous discussion and Lemma \ref{tod2}, SBP$^*$ is the process where each point has a Poisson$(c_1)$ distributed number of offspring and these offspring are distributed bivariate normally with zero mean and the identity covariance matrix. The same process started at $x$ rather than the origin is denoted by SBP$^*_x$.

\begin{lemma}\label{proj}
Given $\epsilon>0$ and $c_1$ sufficiently large, we can find a positive integer $N_0$ and a positive number $M$, such that for all $x \in S_{i,j}=[M(i-\frac{1}{2}),M(i+\frac{1}{2})] \times [M(j-\frac{1}{2}),M(j+\frac{1}{2})]$, the probability that the $N_0$th generation of SBP$^*_x$ contains at least one point in $S_{i+1,j-1}$ and at least one point in $S_{i+1,j+1}$ exceeds $1-\epsilon$.
\end{lemma}

This result is a small perturbation of a result from \cite{hag} and is stated without proof. The process considered in \cite{hag} is a branching random walk, so it is (slightly different) than the process considered here. However, the two processes are very similar and the key difference between the branching random walk and our SBP$^*$ is that the branching random walk can never go extinct. The role of $c_1$ in Lemma \ref{proj} is to make the probability of survival for $N_0$ generations sufficiently large, so that the result of the lemma applies to SBP$^*$ also. We are now ready to prove the main theorem of this section.\\

\noindent{\bf Proof of Theorem \ref{muchwork}: }Fix $\alpha>1$ and consider defining the SBP for a fixed dimension by choosing $\delta_1>0$ as follows. First set $\delta_2>0$ such that $1-\delta_2 > \alpha^{-1}$, and then choose $\delta_1$ such that $1+\delta_1 < (1-\delta_2)\alpha$. Consider running the SBP for $n$ generations and that all points are at least distance $1+\delta_1$ from all others points that aren't its parent. Then, for all $x \in$ SBP, the only points in $S_{1+\delta_1}(x)$ are the offspring of $x$ coming from the homogeneous Poisson process on  $S_{1+\delta_1}(x)$ used to generate them. Thus, in this case, the SBP creates a homogeneous Poisson process on the space $\bigcup S_{1+\delta_1}(x)$, where the union is over all points $x$ in the first $n-1$ generations of the SBP. Furthermore, consider that for all such $x$, the offspring  are born at least distance $1-\delta_2$ away, i.e., all offspring of $x$ appear on the annulus $S_{1+\delta_1}(x)-S_{1-\delta_2}(x)$. This implies that $d_1(x)>1-\delta_2$ for all points $x$ in SBP. Therefore, $r(x)>(1-\delta_2)\alpha>1+\delta_1$, for all $x$ in SBP.

We generate an object that is dominated by the cluster of the site percolation process by running the following algorithm. The algorithm consists of the steps $(0,0),(1,-1),(1,1),(2,-2),\dots,$ with step $(i,j)$ only being carried out if at least one of the steps $(i-1,j-1)$ and $(i-1,j+1)$ is successful. Step $(i,j)$ consists of a SBP started from a point in $S_{i,j}$ and is called {\em successful} unless one of the following errors occur.

\renewcommand{\labelenumi}

\begin{enumerate}
\item{(a)}
The spatial branching process fails to reach the two neighbouring boxes in the projected space.
\item{(b)}
An individual is born within distance $1-\delta_2$ of its parent.
\item{(c)}
An individual in the projected SBP is further than $R_0$ from the origin (of the branching process). The choice of $R_0$ is given below.
\item{(d)}
An individual is born within distance $1+\delta_1$ of an already generated individual that is not its parent.
\end{enumerate}

A step is stopped if any of the above errors occur. To avoid ambiguity, step $(i,j)$ is started from the point in generation $N_0$ of step $(i-1,j-1)$ (if successful) or step $(i-1,j+1)$ (if successful) whose projection is closest to $(Mi,Mj)$. We let ${\cal F}_{i,j}$ denote the $\sigma$-algebra generated by the indicator functions of the successes of steps $(0,0),(1,-1),(1,1),\dots,(i,j-2)$ of the algorithm.

This algorithm generates a point process consisting of all the points contained in all the steps. This point process is then thinned by removing any points that are born in a previously explored part of space. Note that these points are a subset of those that cause type (d) errors and can only occur in unsuccessful steps of the algorithm. This thinning procedure ensures that the point process is a homogeneous Poisson process on a random subset of $\mathbb{R}^d$. Furthermore, we claim that if the probability of a successful step is sufficiently large, this algorithm will with positive probability generate a subset of the GN$_1(d,\alpha)$ model that contains an infinite cluster.

It should be clear from the above discussion that, taken on its own, a successful step generates a homogeneous Poisson process on a random subset of $\mathbb{R}^d$. Furthermore, by the law of GN$_1(d,\alpha)$, the origin of the step is contained in the same cluster as two points that can be used to start two subsequent steps. This still holds after the point process is thinned (since the thinning can only affect points in unsuccessful steps). Furthermore, subsequent steps can only interfere with an earlier successful step with precisely those points that are removed with the thinning. Thus,  Theorem \ref{muchwork} follows from the following claim.
For all sufficiently large $d$, and for all $(i,j) \in {\cal L}$,
$$\mathbb{P}(\text{step }(i,j)\text{ is successful}|{\cal F}_{i,j})=p>p_c.$$

The approach is to show that we can simultaneously make the probabilities of each type of error arbitrarily small by suitable parameter choices. The reason for introducing error (c) is to help bounding error (d) independently of the history of the algorithm. Let $\gamma$ be such that $1-4\gamma>p_c$. We proceed by bounding the probability of each type of error from above by $\gamma$.

Lemma \ref{proj} shows that we can make the probability of error (a) arbitrarily small, for some $c_1$, $N_0$ and for all $d$ sufficiently large. We choose $c_1$ and $N_0$ such that $\mathbb{P}(\text{type (a) error})<\gamma$ for all suitably high dimensions. We also fix $1-\delta_2 > \alpha^{-1}$ and $c_2>c_1$. Recall from the prior discussion about the spatial branching process that for $c_1$ and $c_2$ fixed there exists a $\delta_1(d)$ that defines the required branching process. Therefore, we have also fixed $\delta_1(d)$. Furthermore, since $\delta_1(d) \to 0$ as $d \to \infty$, $1+\delta_1(d)<(1-\delta_2)\alpha$ in high dimensions, as required.

Recall that for all $\epsilon>0,  |S_{1-\epsilon}|/|S_1| \to 0 \text{ as } d \to \infty$. Since the maximum number of points in a step of the algorithm is now bounded by $c_2^{N_0}$ and $\delta_2$ is fixed, we make $\mathbb{P}(\text{type (b) error})$ smaller than $\gamma$ by taking the dimension high enough. Next we choose $R_0$ such that the probability that all individuals of SBP$^*$ are within distance $R_0$ of the origin is at least $1-\gamma$. This implies that in high enough dimensions, the probability of a type (c) error is also less than $\gamma$.

It now only remains to bound the probability of a type (d) error. This is the only type of error that depends on the history of the process and thus bounding this error is more involved. We begin by considering a type (d) error in step $(i,j)$ caused by an individual from step $(\hat{i},\hat{j})$ such that $M\sqrt{(i-\hat{i})^2+(j-\hat{j})^2}<T_0$, where $T_0$ is a constant that will be chosen later on. This is straightforward for any $T_0$ since the number of steps $(\hat i, \hat j)$ to consider is bounded, meaning that the total number of Poisson points is also bounded. Thus, by Lemma \ref{noverlap}, taking the dimension high enough makes the probability of this type of error less than $\gamma/2$.

To finish the proof, we consider a type (d) error due to a previous step $(\hat{i},\hat{j})$ satisfying
\[\lfloor M\sqrt{(i-\hat{i})^2+(j-\hat{j})^2}\rfloor =Q>T_0.\] We define the {\em volume} of the step $(\hat{i},\hat{j})$ to be $|\cup S_{1+\delta_1}(x)|$, where the union is over all points $x$ contained in the step. Note that the volume of step $(\hat{i},\hat{j})$ is uniformly bounded and that the step $(i,j)$ can only fail because of a type (d) error with step $(\hat{i},\hat{j})$ if one of its points falls within the volume of $(\hat{i},\hat{j})$. Note that if an error of type (c) occurs, we stop running the step of the algorithm. So any point that is born a projected distance greater than $R_0$ from the origin of the step has no offspring and is not scanned around. Thus, any point in step $(\hat{i},\hat{j})$ that we scan around must be at least distance $Q-2R_0$ from any point in step $(i,j)$. Thus from Lemma \ref{tod2} and the exponential decay of the normal distribution, for $d$ large enough, the fraction of the projected volume that falls into the circle  with radius $R_0$ centred at $(Mi,Mj)$  is less than $|S_{1+\delta_1}|/Q^3$
for all large $Q$.

Since the number of points $(\hat{i},\hat{j})$ such that $\lfloor M\sqrt{(i-\hat{i})^2+(j-\hat{j})^2}\rfloor =Q$ is bounded by a constant times $Q$, and the series $\sum_{q>T_0}q^{-2}$ converges, we can make the total of such volume from all points small by choosing $T_0$ sufficiently large.  In particular, we can make this volume so small that the probability of each individual in step $(i,j)$ coming within distance $1+\delta_1$ of another point already generated (other than its parent) is less than $\gamma/(2c_2^{N_0})$. Since there are at most $c_2^{N_0}$ points in step $(i,j)$, we get the desired bound.

$\hfill \Box$

\section{Results  for GN$_k(d,\alpha)$ for  $d \ge 2$}\label{subs4}

In low dimensions, the difference between GN$_1$ and GN$_k$ with $k>1$ can be quite pronounced. However, in high dimensions this is no longer the case, since there the difference between $d_k(x)$ and $d_1(x)$ disappears. This leads to some interesting behaviour.

First, we discuss an immediate consequence of Theorem \ref{muchwork}. Recall from Section \ref{sec:2} that for all $\alpha<1$, the GN$_2(d,\alpha)$ is dominated by the GN$_1(d,1)$ model and that there is no percolation in the GN$_1(d,1)$ model for all $d$. Hence, $\alpha_c^{(2)}(d)\geq 1$ for all $d$. Since $\alpha_c^{(2)}\leq \alpha_c$, we have the following corollary to Theorem \ref{muchwork}.

\begin{corollary}
We have $\alpha_c^{(2)}(d) \to 1$  as  $d \to \infty$.
\end{corollary}

In fact, Theorem \ref{muchwork} is expected to hold for the GN$_k(d,\alpha)$ model as well, for all $k$:

\begin{conjecture}\label{conj1}
For all $k \ge 1$, $\alpha_c^{(k)}(d) \to 1 \text{ as } d \to \infty$.
\end{conjecture}

The next result shows that for $d \ge 2$, any $\alpha$ is sufficient for percolation so long as $k$ is suitably large.

\begin{theorem}
For $\alpha>0$ and $d\geq 2$ fixed, there exists a $k$ such that the model with $\alpha_k=\alpha$
percolates on $\mathbb{R}^d$.
\label{smalla}
\end{theorem}

Before giving the proof, we explain the idea of the proof. We divide $\mathbb{R}^2$ into unit squares to create a version of $\mathbb{Z}^2$. A square is called {\em good} if a certain property, depending only on the Poisson process $X$ within the square, holds. In this way we obtain independent site percolation on $\mathbb{Z}^2$, which percolates if the probability that a square is good is greater than $p_c$. It will then be demonstrated that if we have percolation of these good squares, then the underlying continuum percolation model percolates. Finally, we will demonstrate that we can do this in such a way that the probability of a good square is greater than $p_c$.\\ \\
{\bf Proof: } Let $\alpha>0$.
Let $p_c(d)$ denote the critical value of site percolation on $\mathbb{Z}^d$ and recall that $p_c(d)<1$ for all $d \ge 2$ \cite{grim}. The following argument considers the case when $d=2$ and relies on $p_c = p_c(2)<1$. The generalisation to higher dimensions is straightforward and is therefore omitted.

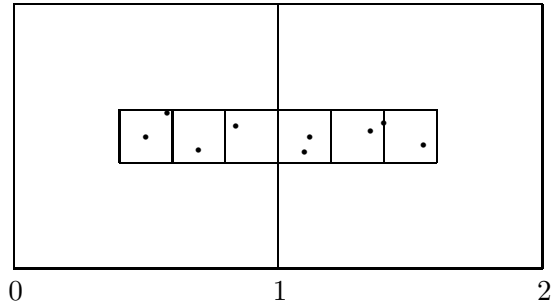
\begin{figure}[htb]
\begin{center}
\begin{picture}(150,150)(0,0)

\put(0, 0){\line(1,0){200}}
\put(0, 100){\line(1,0){200}}
\put(0, 0){\line(0,1){100}}
\put(100, 0){\line(0,1){100}}
\put(200, 0){\line(0,1){100}}
\put(40,40){\line(1,0){120}}
\put(40,60){\line(1,0){120}}
\put(40,40){\line(0,1){20}}
\put(60,40){\line(0,1){20}}
\put(80,40){\line(0,1){20}}
\put(100,40){\line(0,1){20}}
\put(120,40){\line(0,1){20}}
\put(140,40){\line(0,1){20}}
\put(160,40){\line(0,1){20}}

\put(50,50){\circle*{2}}
\put(70,45){\circle*{2}}
\put(84,54){\circle*{2}}
\put(110,44){\circle*{2}}
\put(112,50){\circle*{2}}
\put(135,52){\circle*{2}}
\put(155,47){\circle*{2}}
\put(58,59){\circle*{2}}
\put(140,55){\circle*{2}}
\put(-2,-12){$0$}
\put(98,-12){$1$}
\put(198,-12){$2$}

\end{picture}
\end{center}
\caption{Two good squares and their subsquares $B_1, B_2, \ldots, B_5$ for $n=5$.}
\label{figgie3}
\end{figure}

We divide each unit square equally into $n^2$ smaller squares, where $n$ is some odd integer. Let $X_n$ denote the number of points in a typical subsquare. We call a square {\em good} if all its subsquares satisfy the following property:
$$1\leq X_n \le \frac{m}{n^2},$$
where $m$ will be chosen later.
Now consider two neighbouring good squares. Without loss of generality, let these be $[0,1]^2$ and $[1,2] \times [0,1]$. Define the subsquares
\[ B_i=\Big[\frac{\frac{n-1}{2}}{n}+\frac{i}{n},\frac{\frac{n-1}{2}}{n}+\frac{i+1}{n}\Big] \times \Big[\frac{\frac{n-1}{2}}{n},\frac{\frac{n-1}{2}}{n}\Big], \qquad 0 \le i \le n,\]
see also Figure \ref{figgie3}.
By assumption,  all the $B_i$ are non-empty, and any point in $B_i$ has not more than $m$ neighbours within distance  $(n-1)/(2n)$.
Hence, for all  $k\geq m+1$,  all  points in a subsquare $B_i$ have a connection range of at least $\alpha(n-1)/(2n)$.

Furthermore, points in neighbouring subsquares are at most distance $\sqrt{5}/n$ apart.  By choosing $n > 1+2\sqrt{5}/\alpha$, we ensure that all points in neighbouring subsquares of a good square are connected. As a consequence, all points in neighbouring good squares are connected.  We can now choose $\lambda$, the density of the Poisson process, to be so large that
$$\mathbb{P}(X_n=0)<\frac{1-p_c}{2n^2}.$$
Furthermore, we choose $m$ such that
$$\mathbb{P}(X_n>\frac{m}{n^2})<\frac{1-p_c}{2n^2}.$$
The above relations imply that the probability that a square is good is at least $p_c$. Since squares are good  independently of the state of the other squares, and points in good squares are connected, the good squares dominate independent site percolation on $\mathbb{Z}^2$ with parameter at least $p_c$. Hence, the model percolates.
 $\hfill \Box$\\

Note that Theorem \ref{prop1}, Conjecture \ref{conj1} and Theorem \ref{smalla} would together imply that the critical value is not monotone in the dimension, $d$. By Theorem \ref{smalla}, there exists a $k$ such that $\alpha_c^{(k)}(2) < 0.5$. Theorem \ref{prop1} shows that  $\alpha_c^{(k)}(1)=\infty$ for all $k$. However, Conjecture \ref{conj1} would give that $\alpha_c^{(k)}(d) \to 1$ as $d \to \infty$ for all $k$. In particular, for every $k$ there exists a $d>2$ such that $\alpha_c^{(k)}(d)>0.5$.
Thus, $$\alpha_c^{(k)}(1)>\alpha_c^{(k)}(2)< \alpha_c^{(k)}(d).$$

\section{Results  for GN$(d,\underline{\alpha})$ for $d \ge 2$ }\label{subs5}

We now look at the general model, by removing the restriction that $\underline \alpha$ can have only one non-zero component.

\begin{theorem}
Let $x$ be a point of the Poisson process on $\mathbb{R}^d$. Then $r(x)=\infty$ a.s.\ if and only if  $\sum_{i=1}^\infty i^{1/d}\alpha_i=\infty$.
\end{theorem}
{\bf Proof:} First choose $0<\varepsilon<1$.
Let $(\Gamma_i)$ be a sequence of  disjoint $d$-dimensional annuli centred at $x$ such that for all $k$
the volume of $\Gamma_k$ is $1-\varepsilon$  and $\cup_{i=1}^k \Gamma_i$ is a ball.
Since the volume of a ball of radius $r$ in $\mathbb{R}^d$ is $c(d)r^d$ with $c(d)=\pi^{d/2}/(\Gamma(d/2+1))$,
 the outer radius $\gamma_k$ of $\Gamma_k$ satisfies $c(d)\gamma_k^d=k(1-\varepsilon)$, i.e., $\gamma_k=(\frac{(1-\varepsilon)k}{c(d)})^{1/d}$.

Now define $Y_k$ to be the number of points in $\Gamma_k$, for every $k$, and let $X_k=1-Y_k$. Then the random walk  $S_n=X_1+\cdots +X_n$ has drift $EX_k=1-EY_k=1-$volume$(\Gamma_k)=1-(1-\varepsilon)=\varepsilon>0.$ Since such a  random walk is transient and converges to $\infty$, there exists an a.s.\ finite  random index $N$ such that $S_n>0$ for all $n\geq N$.

If $S_n>0$, then $Y_1+\cdots+ Y_n <n$, and the $n$th nearest neighbour of $x$ is further away than the outer radius $\gamma_n$ of the $n$th annulus, i.e., $d_n(x)\geq \gamma_n= (\frac{(1-\varepsilon)n}{c(d)})^{1/d}$.
Hence, we may write
\[ r(x)=\sum_i \alpha_i d_i(x)\geq  \sum_{i=1}^{N-1} \alpha_i d_i(x)+ \Big(\frac{1-\varepsilon}{c(d)}\Big)^{1/d}\sum_{i=N}^{\infty} \alpha_i i^{1/d}.
\]
So, if $\sum_i \alpha_i i^{1/d}=\infty$, then $r(x)=\infty$ a.s.
Analogously, by considering annuli with volume $1+\varepsilon$, we can find the upper bound
\[ r(x)\leq \sum_{i=1}^{M-1} \alpha_i d_i(x)+ \Big(\frac{1+\varepsilon}{c(d)}\Big)^{1/d}\sum_{i=M}^{\infty} \alpha_i i^{1/d},\]
where $M$ is a.s.\ finite.
Hence, $r(x)\leq \infty$ a.s.\ if $\sum_i \alpha_i i^{1/d}<\infty$. This completes the proof.\hfill $\Box$\\

Since the GN$_k$ model is non-trivial when $d \ge 2$, the behaviour of the model (i.e., percolation or otherwise) depends on both the tail of $\underline{\alpha}$ and the individual $\alpha_i$'s themselves. For this reason, it is clear that there can be no analogy to Conjecture \ref{1d} in the one dimensional case.

\section{Concluding Remarks}

This paper introduces a new continuum percolation model and proves a number of results about this model. In particular, the non-triviality of the model and some features of the critical value have been discussed. Obviously, there are still many open problems relating to the model, some of which are listed below.

\begin{itemize}

\item
The GN$(1,\underline{\alpha})$ model does not percolate if the series $\sum_{i=1}^{\infty} i \alpha_i$ converges?

\item
 If $|\underline{\alpha}| = \sum_{i=1}^{\infty}  \alpha_i =\sum_{i=1}^k \alpha_i< 1$ for some $k$, then does there exist a $d_0$ such that for all $d>d_0$ the GN($d,\underline{\alpha})$ model does not percolate? This is a generalisation of Conjecture \ref{conj1}. A consequence of this result would be the non-monotonicity of the critical value (when varying as a function of the dimension).

\item
It might be possible to give a fuller description of the GN$(d,\underline{\alpha})$ when $|\alpha|=1$. For example, when $d$ is suitably large, the GN$_2(d,1)$ model percolates, but the GN$_1(d,1)$ model does not. How does the GN$(d,\underline{\alpha})$ model behave when $\alpha_1 + \alpha_2 = 1$ and $\alpha_i=0$ for all $i\geq 3$.

\item
This paper has briefly touched upon the related hard sphere models. The author's of this paper would like to again highlight the interesting open problem as to whether percolating hard-sphere models exist in low dimensions.

\end{itemize}

Finally, it would also be interesting to look at other properties of the random graphs created by this model than just percolation.\\

\bibliographystyle{plain}

\end{document}